# Structure of Leavitt path algebras of polynomial growth

Adel Alahmedi[a], Hamed Alsulami[a], Surender Jain[a,b], and Efim I. Zelmanov[a,c,1]

[a]Department of Mathematics, King Abdulaziz University, Jeddah 22254, Kingdom of Saudi Arabia; [b]Department of Mathematics, Ohio University, Athens, OH 45701; and [c]Department of Mathematics, University of California, San Diego, La Jolla, CA 92093

We determine the structure of Leavitt path algebras of polynomial growth and discuss their automorphisms and involutions.

Gelfand–Kirillov dimension | Toeplitz algebra

Following the works in refs. 1–4, Abrams and Aranda Pino (5) and Ara et al. (6) introduced Leavitt path algebras of directed graphs as algebraic analogs of C* algebras of Cuntz and Krieger. This construction provided a rich supply of finitely presented algebras having interesting and extreme properties.

Let $\Gamma = (V, E)$ be a finite directed graph with the set of vertices $V$ and the set of edges $E$. For an edge $e \in E$, we let $s(e)$ and $r(e) \in V$ denote its source and range, respectively. A vertex $v$ for which $s^{-1}(v)$ is empty is called a *sink*. A *path* $p = e_1 \ldots e_n$ in a graph $\Gamma$ is a sequence of edges $e_1 \ldots e_n$ such that $r(e_i) = s(e_{i+1})$, $i = 1, 2, \ldots, n-1$. In this case we say that the path $p$ starts at the vertex $s(e_1)$ and ends at the vertex $r(e_n)$. If $s(e_1) = r(e_n)$ then the path is closed. If $p = e_1 \ldots e_n$ is a closed path and the vertices $s(e_1), \ldots, s(e_n)$ are distinct, then the subgraph $(\{s(e_1), \ldots, s(e_n)\}, \{e_1, \ldots, e_n\})$ of the graph $\Gamma$ is called a *cycle*.

Let $\Gamma$ be a finite graph and let $F$ be a field. The Leavitt path $F$ algebra $L(\Gamma)$ is the $F$ algebra presented by the set of generators $\{v | v \in V\}$ $\{e, e^* | e \in E\}$ and the set of relations (i) $v_i v_j = \delta_{v_i, v_j} v_i$ for all $v_i, v_j \in V$; (ii) $s(e)e = er(e) = e$, $r(e)e^* = e^* s(e) = e^*$ for all $e \in E$; (iii) $e^* f = \delta_{e,f} r(e)$ for all $e, f \in E$; and (iv) $v = \sum_{s(e)=v} ee^*$ for an arbitrary vertex $v \in V \setminus \{\text{sinks}\}$. The mapping that sends $v$ to $v$, $v \in V$, $e$ to $e^*$, and $e^*$ to $e$, $e \in E$, extends to an involution of the algebra $L(\Gamma)$. If $p = e_1 \ldots e_n$ is a path, then $p^* = e_n^* \ldots e_1^*$.

In ref. 7 we showed that the algebra $L(\Gamma)$ has polynomial growth if and only if no two cycles of $\Gamma$ intersect. Let $N = \{1, 2, \ldots\}$, and let $n \in N$. For an algebra $R$, let $M_n(R)$ denote the algebra of $n \times n$ matrices over $R$ and let $M_\infty(R)$ denote the algebra of infinite $N \times N$ finitary matrices over $R$, that is, infinite $N \times N$ matrices with only finitely many nonzero entries.

**Theorem 1.** *Let $L(\Gamma)$ be a Leavitt path algebra of polynomial growth. Then $L(\Gamma)$ has a finite chain of ideals, $(0) \leq I_0 < I_1 < \ldots < I_s = L(\Gamma)$, such that $I_0$ is a finite sum of matrix algebras and infinite finitary matrix algebras over $F$ and each factor $I_{i+1}/I_i$, $i \geq 1$, is a finite sum of matrix algebras and finitary matrix algebras over the Laurent polynomial algebra $F[t^{-1}, t]$. The ideals $I_i$ are invariant under $\text{Aut}(L(\Gamma))$.*

**Remark 1:** We will show that $I_0$ is the locally finite radical of $L(\Gamma)$ (8).

In the rest of the paper we study the algebraic Toeplitz algebra $L(\Gamma_1)$, $\Gamma_1 = \circlearrowright \cdot$ (9) as the simplest nontrivial example of a Leavitt path algebra of polynomial growth. As shown in ref. 9 (it follows also from Theorem 1 above) the locally finite radical $I_0$ of $L(\Gamma_1)$ is $M_\infty(F)$ and $L(\Gamma_1)/M_\infty(F) \cong F[t^{-1}, t]$.

**Theorem 2.** *The short exact sequence $(0) \to M_\infty(F) \to L(\Gamma_1) \to F[t^{-1}, t] \to (0)$ does not split.*

The significance of Theorem 2 is that it shows that the extensions in Theorem 1, generally speaking, do not split.

We describe automorphisms and involutions of the algebraic Toeplitz algebra $L(\Gamma_1)$. Description of involutions is related to the question of whether isomorphic Leavitt path algebras are isomorphic as involutive algebras (10).

**Theorem 3.** $\text{Aut}(L(\Gamma_1)) \cong F^* \ltimes GL_\infty(F)$, *a semidirect product of the multiplicative group $F^*$ of the field $F$ with the general linear finitary group $GL_\infty(F)$. If $F^2 = F$ then the only involution on $L(\Gamma_1)$ (up to isomorphism) is the standard involution $*$.*

In what follows we will assume that the finite graph $\Gamma$ does not have distinct intersecting cycles, which guarantees that $L(\Gamma)$ has polynomial growth. For an arbitrary path $p$, the element $pp^*$ is an idempotent. Consider the family of idempotents $\tilde{\mathcal{E}} = \{pp^* | p \text{ is a path}\}$.

**Remark 2:** We view vertices as paths of length 0.

For two idempotents $e = pp^*$, $f = qq^* \in \tilde{\mathcal{E}}$, if neither $p$ nor $q$ is an initial subpath of the other, then $e$ and $f$ are orthogonal. If $p = qp'$ then $ef = fe = e$.

Consider the set of vertices $V_0 = \{v \in V | \text{ no path starting at } v \text{ finishes at a cycle}\}$. The subset $V_0$ is hereditary and saturated (5). Hence, the ideal $I_0 = \text{id}_{L(\Gamma)}(V_0)$ is the $F$ span of all products $pq^*$, where $p, q$ are paths, $r(p) = r(q) \in V_0$. Let $\mathcal{E} = \{pp^* \in \tilde{\mathcal{E}} | r(p) \in V_0\} \subseteq \tilde{\mathcal{E}}$. Because $pq^* = (pp^*)(pq^*)(qq^*)$ it follows that $I_0 = \mathcal{E} I_0 \mathcal{E}$. Consider also the set of idempotents $\mathcal{E}_S = \{pp^* \in \mathcal{E} | r(p) \text{ is a sink}\}$. We call idempotents from $\mathcal{E}_S$ minimal. Let $v_1, \ldots, v_l$ be all sinks of $\Gamma$. Let $\mathcal{E}_i = \{pp^* \in \mathcal{E} | r(p) = v_i\}$. Clearly, $\mathcal{E}_S = \mathcal{E}_1 \dot\cup \cdots \dot\cup \mathcal{E}_l$, and $\mathcal{E}_i \mathcal{E}_j = (0)$ if $i \neq j$.

**Lemma 4.** (6).

*i) Every idempotent from $\mathcal{E}$ is a sum of minimal idempotents,*

*ii) if $e \in \mathcal{E}_S$ then $eL(\Gamma)e = Fe$,*

*iii) if $e \in \mathcal{E}_i$, $f \in \mathcal{E}_j$ then $\dim_F eL(\Gamma)f = \delta_{ij}$.*

The set $\mathcal{E}_i$ is infinite if and only if there exists a cycle from which one can get to $v_i$. In that case Lemma 4 implies that $\mathcal{E}_i L(\Gamma) \mathcal{E}_i \cong M_\infty(F)$. Otherwise $\mathcal{E}_i L(\Gamma) \mathcal{E}_i \cong M_k(F)$, where $k$ is the number of paths that end at $v_i$. We thus have proved that $I_0$ is isomorphic to a finite sum of matrix algebras and infinite finitary matrix algebras over $F$.

Recall that an algebra is said to be *locally finite dimensional* if every finitely generated subalgebra of it is finite dimensional. The sum of all locally finite-dimensional ideals of an associative algebra $A$ is a *locally finite-dimensional ideal*, which is called the *locally finite-dimensional radical*, denoted by $\text{Loc}(A)$. For further properties of $\text{Loc}(A)$, see ref. 8.

**Lemma 5.** $I_0 = \text{Loc}(L(\Gamma))$.

The ideal $I_0$ is also the socle of the algebra $L(\Gamma)$ (11).

As shown in ref. 5 $L(\Gamma)/I_0 \cong L(\Gamma')$, where $\Gamma' = (V \setminus V_0, E \setminus r^{-1}(V_0))$; the graph $\Gamma'$ does not have sinks. Without loss of generality consider therefore a finite graph $\Gamma$ such that $GK_{\dim} L(\Gamma) < \infty$ and $\Gamma$ does not have sinks, so $I_0 = \text{Loc}(L(\Gamma)) = (0)$.

Recall that an edge $e$ is called an exit from a cycle $C$ if $s(e)$ lies on $C$, but $e$ is not a part of $C$ (5). A cycle without exits will be referred to as an NE cycle. For an arbitrary vertex $v \in V$ there exists a path that starts at $v$ and ends on a cycle, otherwise $v \in V_0$,

---

which contradicts our assumption. Moreover, because distinct cycles of $\Gamma$ do not intersect and all chains of cycles (7) are finite, it follows that for an arbitrary vertex $v \in V$ there exists a path that starts at $v$ and ends on an *NE* cycle.

Consider the set $V_1 = \{v \in V \mid$ a path that starts at $v$ can end only on an *NE* cycle$\}$. The set $V_1$ is obviously hereditary and saturated. Let $C_1, \ldots, C_l$ be all *NE* cycles of $\Gamma$, $\mathcal{E}'_i = \{pp^* \mid p$ is a path, $r(p) \in C_i\}$. Clearly, $\mathcal{E}'_i \mathcal{E}'_j = (0)$ if $i \neq j$. We define $J = \mathrm{id}_{L(\Gamma)}(V_1)$. Then $J = \mathrm{span}(pq^* \mid r(p) = r(q) \in V_1) = \bigoplus_{i=1}^l \mathcal{E}'_i J \mathcal{E}'_i$.

Consider an *NE* cycle $C_i$ with $d_i$ vertices. In ref. 12 it is shown that the subalgebra $L(C_i) = \mathrm{span}(pq^* \mid p, q$ are both paths on the cycle $C_i)$ is isomorphic to $M_{d_i}(F[t^{-1}, t])$.

**Lemma 6.** *Let $e = pp^*, f = qq^* \in \mathcal{E}'_i$. Then $eL(\Gamma)f = pL(C_i)q^*$.*

If the set $\mathcal{E}'_i$ is infinite, which happens if there exists a cycle $C$ different from $C_i$ and a path $p$ such that $s(p) \in C, r(p) \in C_i$, then by *Lemma 6* $\mathcal{E}'_i L(\Gamma) \mathcal{E}'_i \cong M_\infty(F[t^{-1}, t])$. If $|\mathcal{E}'_i| = k$ then $\mathcal{E}'_i L(\Gamma) \mathcal{E}'_i \cong M_k(F[t^{-1}, t])$.

We have proved that the algebra $J$ is isomorphic to a finite direct sum of matrix algebras and infinite finitary matrix algebras over $F[t^{-1}, t]$. The ideal $I_0$ of the algebra $L(\Gamma)$ has been defined. For $i \geq 1$ define $I_i$ via $I_i/I_{i-1} = J(L(\Gamma)/I_{i-1})$. We have got an ascending chain claimed in *Theorem 1*. The ideal $I_0 = \mathrm{Loc}\, L(\Gamma)$ is invariant under $\mathrm{Aut}\, L(\Gamma)$. To prove that the ideals $I_i, i \geq 1$, are invariant we need to obtain an abstract characterization of the ideal $J$.

**Lemma 7.** *The ideal $J$ is the largest ideal of $L(\Gamma)$ with the property that for an arbitrary element $a \in J$, and an arbitrary finite-dimensional subspace $G$ of $L(\Gamma)$ that generates $L(\Gamma)$, there exists a positive constant $k = k(a, G)$ such that $\dim_F(aG^n a) \leq kn$ for $n \geq 1$.*

**Corollary 8.** *Let $\Gamma_1$ and $\Gamma_2$ be finite graphs, and suppose that $\phi: L(\Gamma_1) \to L(\Gamma_2)$ is an isomorphism and that $L(\Gamma_1)$ has polynomial growth. Then $\phi(J(L(\Gamma_1))) = J(L(\Gamma_2))$.*

*Theorem 1* is proved.

We determined the factors $I_{i+1}/I_i$, but the nature of extensions remains unclear. *Theorem 2* implies that generally speaking they do not split.

The algebra $L(\Gamma_1)$ can be presented by generators and relators as $A = \langle x, y \mid xy = 1 \rangle$ (13). Let us fix the notation. The element $e = yx$ is an idempotent. We have

$$I = A(1-e)A = \sum_{i,j=1}^\infty Fe_{ij}, e_{ij} = y^{j-1}(1-e)x^{j-1}, e_{ij}e_{pq} = \delta_{jp}e_{iq}, x(1-e)$$
$$= (1-e)y = 0;$$
$$I \cong M_\infty(F), A/I \cong F[t^{-1}, t].$$

Suppose that the extension splits, that is, the algebra $A$ contains a subalgebra $B$, which is isomorphic to $F[t^{-1}, t]$, $A = B + I$. Let $x = b_1 + \sum_{i,j} \alpha_{ij} e_{ij}$, $1 = b_0 + \sum_{i,j} \beta_{ij} e_{ij}$, $y = b_{-1} + \sum_{i,j} \gamma_{ij} e_{ij}$ where $b_{-1}, b_0, b_1 \in B$; $\alpha_{ij}, \beta_{ij}, \gamma_{ij} \in F$. Consider the finite sets $P(b_1) = \{i \mid \alpha_{ij} \neq 0$ for some $j\}$, $P(b_{-1}) = \{j \mid \gamma_{ij} \neq 0$ for some $i\}$ and one-sided ideal $\rho = \sum_{i \in P(b_1), r \geq 1} Fe_{ir} \triangleleft I$, $\sigma = \sum_{j \in P(b_{-1}), r \geq 1} Fe_{rj} \triangleleft I$.

**Lemma 9.** *For an arbitrary element $a \in I$ there exists $N(a) \geq 1$ such that $b_1^k a \in \rho$, $ab_{-1}^k \in \sigma$, for any $k \geq N(a)$.*

*Proof*: Choose an element $a = \sum_{i,j} \xi_{ij} e_{ij}, 0 \neq \xi_{ij} \in F$. We define $|a| = \max\{i \notin P(b_1) \mid \xi_{ij} \neq 0\}$. Because $xe_{ij} = e_{i-1,j}$ for $i > 1$, $xe_{1j} = 0$, it follows that $|b_1 a| < |a|$ or $a \in \rho$. This implies the first inclusion. The second inclusion is proved in the same way. This completes the proof of the Lemma. ∎

**Lemma 10.** *For an arbitrary element $a \in I$, we have $\dim_F aB < \infty$.*

*Proof*: Let $a \in I$. For arbitrary integers $p, q \geq N(a)$ we have $b_1^p a b_{-1}^q \in \rho \cap \sigma$. Notice that $\dim_F \rho \cap \sigma < \infty$. Fix $p \geq N(a)$. It follows from the above that there exists a nonzero polynomial $f(t) \in F[t]$ such that $b_1^p a f(b_{-1}) = 0$. Because every nonzero ideal of the algebra $F[t^{-1}, t]$ is of finite codimension, we conclude that $\dim_F b_0 aB < \infty$. The element $b_0$ is the identity of the algebra $B$.

Notice that $BI \neq (0)$. Otherwise $IB$ is a nilpotent left ideal of the algebra $I$, which implies that $IB = (0)$, $A = B \oplus I$ is a direct sum. However, the algebra $F[t^{-1}, t] \oplus M_\infty(F)$ is not finitely generated, a contradiction. In view of the simplicity of the algebra $I$, the subset $b_0 I$ generates $I$ as an ideal. This completes the proof of the Lemma.

**Lemma 11.** *For an arbitrary element $a \in I$, we have $\dim_F aA < \infty$.*

*Proof*: Let $a_1$ denote the sum $a_1 = \sum_{i,j} \alpha_{ij} e_{ij}$, and let $a_{-1}$ denote the sum $a_{-1} = \sum_{i,j} \gamma_{ij} e_{ij}, x = b_1 + a_1$, and $y = b_{-1} + a_{-1}$. Let $a \in I$ and let $d = \max(\dim_F aB, \dim_F a_{-1}B, \dim_F a_1 B)$. We claim that for each element $u \in \{a, a_{-1}, a_1\}$ and for an arbitrary product $b$ of elements $b_1, b_{-1}$ of length $d+1$ we have $ub = \sum_k \xi_k ub^{(k)}$, where $\xi_k \in F, b^{(k)}$ are products of lements $b_1, b_{-1}$ of length $\leq d$. Indeed, consider the ascending chain of subspaces $Fu \subseteq uB^{(1)} \subseteq \cdots \subseteq uB^{(d+1)}$, where $B^{(k)}$ is the $F$ span of all products of elements $b_1, b_{-1}$ of length $\leq k$. Because $\dim_F uB \leq d$ we cannot have a strict inclusion at every step. Hence $uB^{(d)} = uB^{(d+1)}$, as claimed. Every product of elements $x, y$ is a linear combination of products of $b_1, b_{-1}, a_1, a_{-1}$. Let $w$ be a product of elements $b_1, b_{-1}, a_1, a_{-1}$. Then $aw$ can be represented as $aw = v_1 w_1 v_2 \ldots v_s w_s$, where $v_i$ are products of $a, a_{-1}, a_1$; $w_i$ are products of $b_1, b_{-1}$. Because of the presence of the element $a$ at the left end the word $v_1$ is not empty. The claim above implies that the words $w_1, \ldots, w_s$ can be assumed to have lengths $\leq d$. Now $aw$ lies in the subalgebra of $M_\infty(F)$ generated by $ab, a_{-1}b, a_1 b$, where elements $b$ are products in $b_1, b_{-1}$ of lengths $\leq d$. This subalgebra is finitely generated, hence finite dimensional. This completes the proof of the Lemma. ∎

It is well known that the set $\{x^i y^j : i, j \geq 0\}$ is a basis of $A$. Hence the elements $(1-yx)y^i, i \geq 0$, are linearly independent, $\dim_F(1-e)A = \infty$, a contradiction. *Theorem 2* is proved.

Now our aim is description of automorphisms and involutions of the algebra $L(\Gamma_1), \Gamma_1 = \bigcirc$.

Consider the countably infinite-dimensional vector space $V = \sum_{i=1}^\infty Fe_i$. Let $E$ be the algebra of all linear transformations of $V$. Because the basis $\{e_i, i \geq 1\}$ has been fixed we can identify $E$ with the algebra of $N \times N$ matrices having only finitely many nonzero entries in each column. Consider also the subalgebra $E_0$ of $E$ which consists of $N \times N$ matrices having finitely many nonzero entries in each row and in each column. As above, $M_\infty(F)$ is the algebra of finitary (having finitely many nonzero entries) $N \times N$ matrices. It is easy to see that $M_\infty(F)$ is an ideal in $E_0$ and a left ideal in $E$.

As follows from *Theorem 1*, the ideal $I_0 = \mathrm{id}_{L(\Gamma_1)}(v_2)$ is isomorphic to $M_\infty(F)$. Extending this isomorphism we can embed $L(\Gamma_1)$ into the algebra $E_0$, the cycle $c$ and its conjugate $c^*$ are identified with the matrices $c = \sum_{i=1}^\infty e_{i+1,i}, c^* = \sum_{i=1}^\infty e_{i,i+1}$, respectively, $e_{ij}$ are matrix units, $L(\Gamma_1) = \,<c, c^*, M_\infty(F)>$.

**Theorem 12.** *(Jacobson, ref. 14). For an arbitrary automorphism $\varphi$ of $M_\infty(F)$ there exists an invertible element $T \in E$ such that $\varphi(a) = T^{-1} aT$ for any $a \in M_\infty(F)$.*

**Lemma 13.** *An automorphism of $L(\Gamma_1)$ induces an automorphism of the type $t \to \alpha t, 0 \neq \alpha \in F, L(\Gamma_1)/I_0 \cong F[t^{-1}, t]$.*

*Proof*: If the assertion is not true then there exists an automorphism $\varphi$ of $L(\Gamma_1)$ whose image in $\mathrm{Aut}\, F[t^{-1}, t]$ maps $t$ to $t^{-1}$. By Jacobson's theorem there exists an invertible element $T \in E$ such that $T^{-1} aT = \varphi(a)$ for all $a \in L(\Gamma_1)$. In particular, $T^{-1} cT = c^* + a, a \in M_\infty(F)$. Hence, $cT = Tc^* + Ta, Ta \in M_\infty(F)$. This implies that for a sufficiently large $n_0 \geq 1$ we have $(cT)_{ij} = (Tc^*)_{ij}$ provided that $i + j \geq n_0$. Therefore, $T_{i-1,j} = T_{i,j-1}$. We showed that $T_{ij} = \alpha_{i+j} \in F$ for $i + j \geq n_0$. The $j$th column of the matrix $T$ intersects all diagonals $\{(i,j) \mid i+j = k\}, k \geq j$. Hence if the sequence $\alpha_k, k \geq 1$, contains infinitely many nonzero entries then every column of $T$ contains infinitely many nonzero entries. Hence the matrix $T$ is finitary, a contradiction. This completes the proof of the Lemma. ∎

**Lemma 14.** *If $T \in E$ is invertible and $T^{-1}M_\infty(F)T = M_\infty(F)$ then $T \in E_0$.*

Recall that the group $GL_\infty(F)$ of invertible matrices from $Id + M_\infty(F)$ is called the finitary general linear group (15). It can be realized as the union $GL_\infty(F) = \bigcup_{n \geq 1} GL_n(F)$.

**Lemma 15.** *Let $\varphi \in Aut\, L(\Gamma_1)$, $\varphi|_{L(\Gamma_1)/I_0} = Id$, $T \in E$, $\varphi(a) = T^{-1}aT$ for any $a \in L(\Gamma_1)$. Then $T = \alpha \cdot Id + a$, $0 \neq \alpha \in F$, $a \in M_\infty(F)$.*

*Proof:* By our assumptions $T^{-1}cT = c + a$, $a \in M_\infty(F)$, or, equivalently, $cT = Tc + Ta$, $Ta \in M_\infty(F)$. Hence for a sufficiently large $n_0 \geq 1$ $(cT)_{ij} = (Tc)_{ij}$, $T_{i+1,j} = T_{i,j-1}$ provided that $i + j \geq n_0$ (we assume that $T_{i,0} = 0$). Hence $T$ is an almost Toeplitz matrix, $T = T_0 + \sum_{u \in Z} \alpha_k c^{(k)}$, where $T_0 \in M_\infty(F)$, $c^{(k)} = \sum_{j-i=k} e_{ij}$, $\alpha_k \in F$. The $j$th column intersects all diagonals $\{(i,j)|j-i=k\}$ with $k \leq j$. Hence the set $\{k < 0 | \alpha_k \neq 0\}$ is finite. Similarly, an $i$th row intersects all diagonals $\{(i,j)|j-i=k\}$ with $-k \leq i$. Hence the set $\{k > 0 | \alpha_k \neq 0\}$ is finite as well. Now we have $T = T_0 + \sum_{k=-m}^{n} \alpha_k c^{(k)}$; $\alpha_{-m} \cdot \alpha_n \neq 0$. Because the matrix $\sum_{n=-m}^{n} \alpha_k c^{(k)}$ cannot be strictly upper or lower triangular (otherwise $T$ would not be invertible), we can assume that $m, n \geq 0$. All of the above applies to the matrix $T^{-1}$ as well, $T^{-1} = (T^{-1})_0 + \sum_{s=-p}^{q} \beta_s c^{(s)}$; $\beta_q \cdot \beta_{-p} \neq 0$; $p, q \geq 0$, $(T^{-1})_0$ is a finitary matrix. Now

$$Id = \left(T_0 + \sum_{k=-m}^{n} \alpha_k c^{(k)}\right)\left((T^{-1})_0 + \sum_{s=-p}^{q} \beta_s c^{(s)}\right)$$
$$= T_0 \cdot T^{-1} + (T-T_0)(T^{-1})_0 + \sum \alpha_i \beta_j c^{(i+j)}.$$

Because $T, T^{-1} \in E_0$ it follows that $T_0 \cdot T^{-1} + (T-T_0)(T^{-1})_0 \in M_\infty(F)$. Moreover, the equality above implies that $m = n = p = q = 0$, $T = \alpha_0 \cdot Id + T_0$. This completes the proof of the Lemma, and thus completes the proof of *Theorem 3*. ∎

Consider the embedding $\pi : F^* \to E_0^*$ of the multiplicative group of the field $F$ into the multiplicative group of the algebra $E_0$, $\pi(\alpha) = diag(1, \alpha, \alpha^2, \ldots)$. It is easy to see that $\pi(\alpha)^{-1}L(\Gamma_1)\pi(\alpha) = L(\Gamma_1)$ and $\pi(\alpha)^{-1}c\pi(\alpha) = \alpha c$. Now, *Lemmas 13* and *15* imply that $Aut(L(\Gamma_1)) = \pi(F^*) \ltimes GL_\infty(F)$.

We say that two involutive algebras $(R_1, *_1)$ and $(R_2, *_2)$ are isomorphic if there exists an isomorphism $\varphi : R_1 \to R_2$ of algebras $R^1, R^2$, such that $\varphi(a^{*_1}) = \varphi(a)^{*_2}$ for an arbitrary element $a \in R_1$.

**Lemma 16.** *Let $F^2 = F$. Then the algebra $L(\Gamma_1)$ has only one (up to isomorphism) involution: the standard involution $*$.*

*Proof*: If we view $L(\Gamma_1)$ as a subalgebra of the algebra $E_0$, then the standard involution $*$ becomes the restriction of the transposition $(a_{ij})^t = (a_{ji})$. Let $: L(\Gamma_1) \to L(\Gamma_1)$ be an involution. The composition of the involutions $-$ and $t$ is an automorphism. Hence there exists a matrix $T \in \pi(F^*)GL_\infty(F)$ such that $(\bar{a})^t = T^{-1}aT$ for all elements $a \in L(\Gamma_1), \bar{a} = T^t a^t (T^t)^{-1}$. Applying the involution $-$ twice we get $a = \bar{\bar{a}} = T^t(T^{-1}aT)(T^t)^{-1} = (T^tT^{-1})a(T(T^t)^{-1})$. Because the matrix $(T^t)T^{-1}$ commutes with an arbitrary matrix from $M_\infty(F)$ it follows that $T^tT^{-1} = \alpha \cdot Id$, $\alpha \in F^*$, $T^t = \alpha T$. Now, $T = (T^t)^t = \alpha^2 T$, $\alpha = \pm 1$. All nonzero entries of the matrix $T$ except finitely many lie in the main diagonal. Hence $T$ cannot be skew-symmetric. Hence $T^t = T$. If an arbitrary element from $F$ is a square then there exists a matrix $Q \in \pi(F^*)GL_\infty(F)$ such that $T = Q^tQ$. Now the mapping $a \to Q^{-1}aQ$ is an isomorphism of the involutive algebra $(L(\Gamma_1), t)$ to the involutive algebra $(L(\Gamma_1), *)$. This completes the proof of the Lemma. ∎

**Note Added in Proof.** For a different approach to automorphisms of the Jacobson algebra, see ref. 16.

**ACKNOWLEDGMENTS.** The authors thank G. Abrams and J. Bell for numerous helpful remarks. This research is supported by the Deanship of Scientific Research, King Abdulaziz University. The work is also partially supported by the National Science Foundation.

1. Cuntz J (1977) Simple C*-algebras generated by isometries. *Commun Math Phys* 57(2): 173–185.
2. Cuntz J, Krieger W (1981) A class of C*-algebras and topological Markov chains. *Invent Math* 63(2):25–40.
3. Kumjian A, Pask D, Raeburn I (1998) Cuntz-Krieger algebras of directed graphs. *Pac J Math* 184(1):161–174.
4. Raeburn I (2005) Graph algebras. *CBMS Regional Conference Series in Mathematics* 103 (American Mathematical Society, Providence, RI).
5. Abrams G, Aranda Pino G (2005) The Leavitt path algebra of a graph. *J Algebra* 293(2):319–334.
6. Ara P, Moreno MA, Pardo E (2007) Nonstable K-theory for graph algebras. *Algebr Represent Theory* 10(2):157–178.
7. Alahmadi A, Alsulami H, Jain SK, Zelmanov E (2012) Leavitt path algebras of finite Gelfand-Kirillov dimension. *J Algebra Appl* 11(6).
8. Zhevlakov KA, Shestakov IP (1973) Local finiteness in the sense of Shirshov. *Algebra Log* 12:41–73.
9. Siles Molina M (2008) Algebras of quotients of path algebras. *J Algebra* 319(12): 5265–5278.
10. Abrams G, Tomforde M (2011) Isomorphisms and Morita equivalence of graph algebras. *Trans AMS* 363(7):3733–3767.
11. Aranda Pino G, Martin Barquero D, Martin Gonzalez C, Siles Molina M (2010) Socle theory for Leavitt path algebras of arbitrary graphs. *Rev Mat Iberoam* 26(2):611–638.
12. Abrams G, Aranda Pino G, Perera F, Siles Molina M (2010) Chain conditions for Leavitt path algebras. *Forum Math* 22(1):95–114.
13. Jacobson N (1950) Some remarks on one-sided inverses. *Proc AMS* 1:352–355.
14. Jacobson N (1956) Structure of rings. *AMS Coll Publ* 37 (American Mathematical Society, Providence, RI).
15. Hall JI (2006) Periodic simple groups of finitely linear transformations. *Ann Math* 163(2):445–498.
16. Bavula V (2012) The group of automorphisms of the Jacobian algebra. *J Pure Appl Algebra* 216(3):535–564.